\documentclass[final]{elsarticle}

\usepackage{amssymb,latexsym}
\usepackage{amsmath}

\biboptions{sort&compress}

\journal{: \, J. Comb. Number Theory}

\begin{document}

\newtheorem{teo}{Theorem}
\newproof{prova}{Proof}

\begin{frontmatter}

\title{Another elementary proof of $\: \sum_{n \ge 1}{1/{n^2}} = \pi^2/6\,$ and a recurrence formula for $\,\zeta{(2k)}$}

\author{F. M. S. Lima}

\address{Institute of Physics, University of Bras\'{i}lia, P.O. Box 04455, 70919-970, Bras\'{i}lia-DF, Brazil}

\ead{fabio@fis.unb.br}


\begin{abstract}
In this shortnote, a series expansion technique introduced recently by Dancs and He for generating Euler-type formulae for odd zeta values $\:\zeta{(2 k +1)}$, $\zeta{(s)}$ being the Riemann zeta function and $\,k\,$ a positive integer, is modified in a manner to furnish the even zeta values $\,\zeta{(2k)}$.  As a result, I find an elementary proof of $~\sum_{n=1}^\infty{{1/{n^2}}} = {\,\pi^2/6}\,$, as well as a recurrence formula for $\,\zeta{(2k)}\,$ from which it follows that the ratio ${\, \zeta{(2k)} / \pi^{2k}}\,$ is a rational number, without making use of Euler's formula and Bernoulli numbers.
\newline
\end{abstract}


\end{frontmatter}

\section{Introduction}

For real values of $\,s$, $s>1$, the Riemann zeta function is defined as $\,\zeta(s) := \sum_{n=1}^\infty{{\,1/n^s}}$.\footnote{In this domain, this series converges according to the integral test. For $\,s=1$, one has the harmonic series $\:\sum_{n=1}^\infty{1/n}$, which diverges to infinity.} For $\,s=2k$, $\,k \in \mathbb{Z}$, $k>0$, Euler (1740) did find that~\cite{Euler}  
\begin{equation}
\zeta(2k) = \frac{2^{2k-1} \, \left|B_{2k}\right|}{(2k)!} \, \: \pi^{2k} \, ,
\label{eq:Euler}
\end{equation}
where $\,B_{k}\,$ is the $k$-th Bernoulli number.\footnote{The (rational) numbers $B_{k}$ are the coefficients of ${\,z^k/k!}\,$ in the Taylor series expansion of ${\,z/(e^z-1)}$, $|z| < 2\,\pi$.}
As a consequence, since $\,B_2 = 1/6\,$ one has $\,\zeta(2) = \pi^2/6$, which is the Euler solution to the Basel problem (see Ref.~\cite{Euler0} and references therein).

By noting that the series expansion approach introduced by Dancs and He (2006) on seeking for an Euler-type formula for $\,\zeta{(2k+1)}$, see Ref.~\cite{Dancs}, could be modified in a manner to furnish similar formulas for $\,\zeta(2k)$, here in this note I show that the substitution of $\,\sin{(n \pi)}\,$ by $\,\cos{(n \pi)}\,$ in the Dancs-He initial series in fact yields a series expansion which can be reduced to a finite sum involving only even zeta values.  From the first few terms of this sum, I have found an elementary proof of $\,\zeta{(2)} = {\,\pi^2/6}\,$ and a recurrence formula for $\zeta{(2k)}$.  The proofs are elementary in the sense they do not involve complex analysis, Fourier series, or multiple integrals.\footnote{For non-elementary proofs, see, e.g., Refs.~\cite{Kalman,Apostol} and references therein.}

\section{Elementary evaluation of $\,\zeta{(2)}$}

For any real $\epsilon > 0$ and $u \in [1,1+\epsilon]$, we begin by taking into account the following Taylor series expansion considered by Dancs and He in Ref.~\cite{Dancs}:
\begin{equation}
\frac{2 \, e^t}{e^t + u} = \sum_{m=0}^{\infty}{\phi_m(u) ~ \frac{t^m}{\,m!}} \, ,
\label{eq:phiserie}
\end{equation}
which converges absolutely for $\,|t| < \pi$.

From the generating function for the Euler polynomial $\,E_m(x)$, namely ${\,2 \, e^{x\,t}/(e^t + 1)} = \sum_{m=0}^{\infty}{E_m(x)\, \frac{t^m}{m!}}\,$, it is clear that $\,\phi_m(1) = E_m(1)$, for all nonnegative integer values of $\,m$.  For $u>1$, we have
\begin{equation}
\phi_m(u) = -2 \, \sum_{n=1}^{\infty}{\frac{n^{\,m}}{(-u)^{\,n}}} \, .
\label{eq:phi}
\end{equation}
Let us take this series as our definition of $\,\phi_{-m}(u)$, $m$ being a positive integer. Therefore
\begin{equation}
\phi_{-m}(1) = -2 \, \sum_{n=1}^\infty{\frac{(-1)^n}{n^m}} = -2 \: \zeta^{*}{(m)} = 2\,(1-2^{1-m})\,\zeta{(m)}
\label{eq:phinegative}
\end{equation}
for all integer $\,m>1$.

Now, let
\begin{equation*}
f(u) := \sum_{n=1}^{\infty}{\frac{\left({\,1/u}\right)^n}{n^2}}
\end{equation*}
be an auxiliary function, with $u$ belonging to the same domain as above.  Since $\,\cos{(n \pi)} = (-1)^n$, then $\,f(u)\,$ can be written in the form
\begin{equation*}
f(u) = \sum_{n=1}^{\infty}{(-1)^n \, \frac{\cos(n\,\pi)}{u^n\,n^2}} \, .
\end{equation*}
On expanding $\,\cos{(n \pi)}\,$ in a Taylor series, one has
\begin{equation*}
f(u) = \sum_{n=1}^{\infty}{\left[\frac{(-1)^n}{u^n\,n^2} \cdot \sum_{j=0}^{\infty}{(-1)^j \, \frac{(n \pi)^{2 j}}{(2 j)!}}\right]} = \sum_{j=0}^{\infty}{(-1)^j \frac{\pi^{2 j}}{(2 j)!} \, \sum_{n=1}^{\infty}{(-1)^n \, \frac{n^{2 j}}{u^n\,n^2}}} ,
\end{equation*}
in which the change of sums justifies by Fubini's theorem.  By writing the last series in terms of $\,\phi_m(u)$, one has
\begin{eqnarray}
f(u) = \sum_{j=0}^{\infty}{(-1)^j \, \frac{\pi^{2 j}}{(2 j)!} \, \frac{\phi_{2j-2}(u)}{(-2)}} \, .
\label{eq:marcia}
\end{eqnarray}

This is sufficient for proofing our first result.
\newline

\begin{teo}[Short evaluation of $\,\zeta{(2)}\,$]
\label{teo:z2}
\begin{equation*}
\sum_{n=1}^\infty{\frac{1}{\,n^2}} = \frac{\pi^2}{6} \, .
\end{equation*}
\end{teo}

\begin{prova}
\; By taking the limit as $u \rightarrow 1^{+}$ on both sides of Eq.~\eqref{eq:marcia}, one has
\begin{eqnarray}
\lim_{u \rightarrow 1^{+}} \sum_{n=1}^{\infty}{\frac{1}{u^n \, n^2}} = -\frac12 \, \phi_{-2}(1) \, + \frac12 \, \frac{\pi^2}{2!} \, \phi_0{(1)} -\frac12 \, \sum_{j=2}^{\infty}{(-1)^j \frac{\pi^{2 j}}{(2 j)!} \, \phi_{2j-2}(1)} \, ,
\end{eqnarray}
which, in face of the value of $\phi_{-2}(1)$ stated in Eq.~\eqref{eq:phinegative}, implies that
\begin{equation}
\sum_{n=1}^{\infty}{\frac{1}{n^2}} = -\frac12 \, \left[ 2 \left(1-2^{1-2}\right) \zeta{(2)} \right] \, +\frac{\pi^2}{4} \, E_0{(1)} -\frac12 \, \sum_{j=2}^{\infty}{(-1)^j \frac{\pi^{2 j}}{(2 j)!} \, E_{2j-2}(1)} \, .
\end{equation}
Since $E_0(1)=1$ and $E_m(1)=0$ for all $m>0$, the right-hand side of this equation reduces to $\,-\frac12 \, \zeta{(2)} +{\,\pi^2 / 4}$, which implies that
\begin{equation*}
\zeta{(2)} = -\frac12 \, \zeta{(2)} +\frac{\pi^2}{4} \, ,
\end{equation*}
and then $\: \dfrac{3}{2} \: \zeta{(2)} = \dfrac{\pi^2}{4}\,$.  
\begin{flushright} $\Box$ \end{flushright}
\end{prova}

\section{Recurrence formula for $\zeta{(2 k)}$}

Interestingly, our approach can be easily adapted to treat higher even zeta values by changing the exponent of $\,n\,$ from $\,2\,$ to $\,2 k$. The result is the following recurrence formula for even zeta values.

\begin{teo}[Recurrence for $\,\zeta{(2 k)}\,$]
\label{teo:z2k}
\;  For any positive integer $\,k$,
\begin{equation*}
\left( 4 -\frac{4}{2^{2k}} \right) \zeta{(2k)} =  \sum_{m=1}^{k-1}{ \frac{(-1)^{k-m+1}}{(2k-2m)!} \left( 2 -\frac{4}{2^{2m}} \right) \pi^{2k-2m} \, \zeta{(2m)}} \,-(-1)^k \frac{\pi^{2k}}{(2k)!} \, .
\label{eq:recorrencia}
\end{equation*}
\end{teo}

\begin{prova}
\; We begin by defining $\,f_k{(u)} :=  \sum_{n=1}^{\infty}{\left({\,1/u}\right)^n / n^{2 k}}\,$.  Again, since $\,\cos{(n \pi)} = (-1)^n$, we may write
\begin{eqnarray}
f_k{(u)} &=& \sum_{n=1}^{\infty}{(-1)^n \, \frac{\cos(n\,\pi)}{u^n\,n^{2k}}} = \sum_{n=1}^{\infty}{\frac{(-1)^n}{u^n\,n^{2k}} \, \sum_{j=0}^{\infty}{(-1)^j \, \frac{(n \pi)^{2 j}}{(2 j)!}}} \nonumber \\
&=& \sum_{j=0}^{\infty}{(-1)^j \frac{\pi^{2 j}}{(2 j)!} \, \sum_{n=1}^{\infty}{(-1)^n \, \frac{n^{2 j}}{u^n\,n^{2k}}}}\, .
\end{eqnarray}
On rewriting the last series in terms of $\phi_m(u)$, one finds
\begin{equation*}
f_k{(u)} = \sum_{j=0}^{\infty}{(-1)^j \frac{\pi^{2 j}}{(2 j)!} \, \frac{\phi_{2j-2k}(u)}{(-2)}}  \, -\frac12 \, \sum_{j=0}^{k-1}{(-1)^j \frac{\pi^{2 j}}{(2 j)!} \, \phi_{2j-2k}(u)} \,- \frac12 \, \sum_{j=k}^\infty{(-1)^j \frac{\pi^{2 j}}{(2 j)!} \, \phi_{2j-2k}(u)} \, . 
\end{equation*}
Now, on substituting $\,m = j -k\,$ in the above series, one has
\begin{eqnarray}
f_k{(u)} = -\frac12 \, \sum_{m=-k}^{-1}{(-1)^{m+k} \, \frac{\pi^{2m+2k}}{(2m+2k)!} \, \phi_{2m}(u)} -\frac12 \, \sum_{m=0}^\infty{(-1)^{m+k} \, \frac{\pi^{2m+2k}}{(2m+2k)!} \, \phi_{2m}(u)} \nonumber \\
=  -\frac12 \, (-1)^k \left[ \, \sum_{\widetilde{m}=1}^{k}{ \frac{(-1)^{\widetilde{m}}\,\pi^{2k-2\widetilde{m}}}{(2k-2\widetilde{m})!} \, \phi_{-2\widetilde{m}}(u)} 
+ \sum_{m=0}^\infty{\frac{(-1)^{m} \, \pi^{2m+2k}}{(2m+2k)!} \, \phi_{2m}(u)} \right] \! . \;
\label{eq:2seriesM}
\end{eqnarray}
The limit as $u \rightarrow 1^{+}$, taken on both sides of Eq.~\eqref{eq:2seriesM}, yields
\begin{equation}
\lim_{u \rightarrow 1^{+}} \sum_{n=1}^{\infty}{\frac{1}{u^n \, n^{2k}}} = -\frac12 \, (-1)^k \left[ \, \sum_{m=1}^k{ \frac{(-1)^m \, \pi^{2k-2m}}{(2k-2m)!} \, \phi_{-2m}(1)} 
+\sum_{m=0}^\infty{\frac{(-1)^m \, \pi^{2m+2k}}{(2m+2k)!} \, \phi_{2m}(1)} \right]\!.
\label{eq:aux}
\end{equation}
From Eq.~\eqref{eq:phinegative}, one knows that $\phi_{-2m}(1) = 2 \left( 1 -2^{1-2m} \right) \zeta{(2m)}$. For nonnegative values of $m$, one has $\phi_{2m}(1) = E_{2m}(1) = 0$, the only exception being $\,\phi_0(1) = E_0(1)=1$. This reduces Eq.~\eqref{eq:aux} to
\begin{equation*}
\sum_{n=1}^{\infty}{\frac{1}{n^{2k}}} = -(-1)^{k} \sum_{m=1}^k{\frac{(-1)^m \, \pi^{2k-2m}}{(2k-2m)!} \, \left( 1 -2^{1-2m} \right) \zeta{(2m)}} -(-1)^k\, \frac{\pi^{2k}}{2\,(2k)!} \, .
\end{equation*}
By extracting the last term of the sum and isolating $\,\zeta{(2k)}$, one finds
\begin{equation*}
\left( 2 -\frac{2}{2^{2k}} \right) \, \zeta{(2k)} =  (-1)^{k+1} \sum_{m=1}^{k-1}{\frac{(-1)^m \, \pi^{2k-2m}}{(2k-2m)!} \, \left( 1 -2^{1-2m} \right) \zeta{(2m)}} -(-1)^k\, \frac{\pi^{2k}}{2\,(2k)!} \, .
\end{equation*}
A multiplication by $2$ on both sides yields the desired result. 
\begin{flushright} $\Box$ \end{flushright}
\end{prova}

The first few even zeta values can be readily obtained from this recurrence formula. For $\,k=1$, the sum in the right-hand side is null and one has
\begin{equation*}
3 \,\zeta(2) = -(-1)\,\frac{\pi^2}{2} \, ,
\end{equation*}
which simplifies to $\,\zeta(2) = \pi^2/6$, in agreement to Theorem~\ref{teo:z2}.

For $k=2$, one has
\begin{equation*}
\frac{15}{4} \: \zeta(4) = \frac{\pi^2}{2!} \: \zeta(2) - \frac{\pi^4}{ 4!} \, .
\end{equation*}
By substituting the value of $\,\zeta(2)$, above, and multiplying both sides by $4$, one finds
\begin{equation}
15 \: \zeta(4) = \frac{\pi^4}{3} -\frac{\pi^4}{6} = \frac{\pi^4}{6} \, ,
\label{eq:z4}
\end{equation}
which implies that $\,\zeta(4) = \pi^4 / 90$.

Note that, by writing the recurrence formula in Theorem~\ref{teo:z2k} in the form
\begin{equation}
\left( 1 -\frac{1}{2^{2k}} \right) \, \frac{\zeta{(2k)}}{\pi^{2k}} = \sum_{m=1}^{k-1}{ \, \frac{(-1)^{k-m+1}}{(2k-2m)!} \left( \frac12 -\frac{1}{2^{2m}} \right) \frac{\zeta{(2m)}}{\pi^{2m}}} -\frac{(-1)^k}{4\,(2k)!} \, ,
\end{equation}
it is straightforward to show, by induction on $k$, that the ratio ${\, \zeta{(2k)} / \pi^{2k}}\,$ is a rational number for every positive integer $\,k$, without making use of Euler's formula for $\zeta{(2k)}$, see Eq.~\eqref{eq:Euler}, and Bernoulli numbers. In fact, this was the original motivation that has led the author to study the properties of the Dancs-He series expansions. The proofs developed here could well be modified to cover other special functions of interest in analytic number theory.


\end{document}